\documentclass[12pt,reqno]{amsart} 
\usepackage{amssymb,amsmath}
\usepackage[all]{xy}

\swapnumbers

\DeclareMathOperator{\Coh}{Coh}
\DeclareMathOperator{\Ob}{Ob}
\DeclareMathOperator{\ml}{mod--\Lambda}
\DeclareMathOperator{\Con}{Con}
\DeclareMathOperator{\C}{C}
\DeclareMathOperator{\D}{D}
\DeclareMathOperator{\Hh}{H}
\DeclareMathOperator{\K}{K}
\DeclareMathOperator{\Ll}{L}
\DeclareMathOperator{\T}{T}
\DeclareMathOperator{\gr}{gr}

\DeclareMathOperator{\im}{Im}
\DeclareMathOperator{\Ker}{Ker}

\DeclareMathOperator{\id}{id}
\DeclareMathOperator{\Hom}{Hom}
\DeclareMathOperator{\Ext}{Ext}

\DeclareMathOperator{\rank}{rank}

\newcommand{\Z}{{\mathbb Z}}

\renewcommand{\P}{{\mathbb P}}
\renewcommand{\H}{{\text{H}}}
\renewcommand{\L}{{\text{L}}}

\newcommand{\ka}{{\mathcal A}}
\newcommand{\kb}{{\mathcal B}}
\newcommand{\kc}{{\mathcal C}}

\newcommand{\ke}{{\mathcal E}}
\newcommand{\kf}{{\mathcal F}}
\newcommand{\kg}{{\mathcal G}}
\newcommand{\kh}{{\mathcal H}}
\newcommand{\ki}{{\mathcal I}}

\newcommand{\kl}{{\mathcal L}}

\newcommand{\ko}{{\mathcal O}}
\newcommand{\kp}{{\mathcal P}}

\newcommand{\kv}{{\mathcal V}}

\newcommand{\ky}{{\mathcal Y}}
\newcommand{\kz}{{\mathcal Z}}

\newcommand{\lra}{\longrightarrow}
\newcommand{\lla}{\longleftarrow}

\begin{document}

\title[Horrocks theory and BGG correspondence]{Horrocks theory and the\\
Bernstein-Gel'fand-Gel'fand correspondence }

\author[Coand\u{a}]{I.~Coand\u{a}$^{1}$} \address{
  Institute of Mathematics
  of the Romanian Academy\\
  P.O.~Box 1-764
\newline RO--70700 Bucharest, Romania}
\email{Iustin.Coanda@imar.ro} 
\author[Trautmann]{G.~Trautmann$^{1}$}
\address{
  Universit\"at Kaiserslautern\\
  Fachbereich Mathematik\\
  Erwin-Schr\"odinger-Stra{\ss}e
\newline D-67663 Kaiserslautern}
\email{trm@mathematik.uni-kl.de}

\footnotetext[1]{Research for this paper was partially supported by DFG.}

\begin{abstract}
  We construct an explicit  equivalence between a category of
  complexes over the exterior algebra, which we call HT--complexes, and
  the stable category of vector bundles on the corresponding projective space,
  essentially translating into more fancy terms the results of Trautmann
  \cite{tr2} which, in turn, were influenced by ideas of Horrocks \cite{ho1}, 
  \cite{ho2}.   However, the result expressed by 
  Theorem \ref{stabism} and its corollary, which
  establishes a relation between the Tate resolutions over the exterior
  algebra, described in \cite{efs}, and HT--complexes, might be new,
  although, perhaps, not a surprise to experts.
\end{abstract}

\maketitle
\thispagestyle{empty}
\tableofcontents
\vskip1cm

\section*{Introduction}

The Bernstein--Gel'fand--Gel'fand correspondence states that the derived
category ${\D}^b(\Coh\P_n)$ of coherent sheaves on projective $n$--space is
equivalent to the stable category of the category $\ml$ of finitely
generated graded modules over the exterior algebra $\Lambda = \Lambda^\bullet
V$ of the corresponding vector space $V$ overlying $\P_n$, see \cite{bgg}.
New light was shed onto this correspondence by the paper \cite{efs} of
Eisenbud--Fl{\o}ystad--Schreyer, who studied systematically the so--called
Tate resolutions of sheaves or their graded modules.   These Tate resolutions
are doubly unbounded acyclic complexes of free graded $\Lambda$--modules, see
Section \ref{section5} for a definition.   The main result in \cite{efs} is
that the terms of a Tate resolution of a sheaf are determined by its
cohomology, as well as the linear parts of the differentials.

In \cite{co}, I.~Coand\u{a} showed that a combination of a remark in \cite{bgg}
and the use of Tate resolutions over the exterior algebra leads to quick 
proofs of the main results of both \cite{bgg} and \cite{efs}.

In \cite{tr2}, G.~Trautmann had established a correspondence between complexes 
of extensions of Koszul operators and stable isomorphism classes of vector 
bundles on $\P_n$.   This correspondence leads to applications on the 
structure of stable isomorphism classes of vector bundles analogous to the 
applications of the Beilinson monads.   The complexes in \cite{tr2} resemble 
filtrations studied in \cite{ho2}, and their construction was influenced 
by ideas of G.~Horrocks.

In this paper, we show that the complexes in \cite{tr2} correspond to bounded
complexes $G^\bullet$ of free modules 
over the exterior algebra, of a structure similar to that
of Tate resolutions, called HT--complexes, see \ref{ht}.   It is shown that
the BGG--functor $\L$ defines an equivalence between the homotopy category
$\kh$ of HT--complexes and the stable category of vector bundles over $\P_n$,
see \ref{main}.   In the course of proof, the results of \cite{tr2} are
re--proved.   Furthermore, it is shown that an HT--complex $G^\bullet$ is
determined up to isomorphism by the stable isomorphism class of 
the bundle it is defining.

The other main result of this paper is that an HT--complex $G^\bullet$ of a
bundle is related to the Tate resolution $I^\bullet$ of the bundle by the
formula $G^\bullet = F_{n-1} I^\bullet/F_0 I^\bullet$, see Theorem 
\ref{stabism}.   This reflects the fact that the cohomology $\H^0$ and $\H^n$ 
of bundles is neglected in the consideration of stable isomorphism classes. 

We are grateful to Wolfram Decker for suggesting us that the results from 
\cite{tr2} could be related to the BGG correspondence via the results 
from \cite{efs} about Tate resolutions. We would like, however, to point 
out that, except for the proof of the comparison theorem \ref{stabism}, 
our treatment of the subject is independent of the results from
\cite{bgg} and \cite{efs}.
\vskip1cm

\section{$\Lambda$--modules and associated complexes}

\textbf{Notation.}\label{nota}
We shall use the notations and conventions from \cite{co}.   It will be also 
convenient for us to refer to that paper for the proof of some standard 
facts. 

Let $k$ be a field, $V$ an $(n+1)$--dimensional $k$--vector space, $n \ge 2$,
and let $\P = \P_n(k) = \P V$ denote the projective space of 1--dimensional
subspaces of $V$, such that $\H^0 \ko_{\P V} (1) = V^\ast$.  We shall write
$\ko(d) = \ko_\P (d)$ for the standard invertible sheaf on $\P$ of degree $d$.
$\Lambda = \Lambda(V) = {\oplus}_{d \ge 0} \Lambda^d V$
denotes the
exterior algebra of $V$ and $S = S(V^\ast) = {\oplus}_{d\ge 0} S^d V^\ast$ 
denotes the symmetric algebra of $V^\ast$.

The category of finitely generated graded right $\Lambda$--modules is denoted 
by $\ml$.

The category $\C(\ml)$ of complexes of modules in $\ml$ will be denoted
$\C(\Lambda)$, and the homotopy, resp.\, derived category of $\ml$ by
$\K(\Lambda)$, resp.\, $\D(\Lambda)$.   Similarly, we use the notations
$\C^b(\Lambda)$, $\K^b(\Lambda)$ and $\D^b(\Lambda)$ for the categories of
bounded complexes or $\C^\pm(\Lambda)$, $\K^\pm(\Lambda)$, $\D^\pm(\Lambda)$ 
for the categories of half--bounded complexes.

Similarly we write $\C(\P)$, $\K(\P)$, $\D(\P)$ for the categories 
$\C(Q\Coh(\P))$, $\K(Q\Coh(\P))$, $\D(Q\Coh(\P))$
based on the category of quasi--coherent sheaves over $\P$.
\vskip5mm

\subsection{Free $\Lambda$-modules}\label{freem}
The module $\Lambda$ has the natural positive grading with $\Lambda_i =
\Lambda^iV$, however, the dual module $\Lambda^{\vee}$ is defined by
$\Lambda^{\vee}_i = \Lambda^{-i} V^\ast$.  For an arbitrary $\Lambda$--module
$N$ the dual $N^{\vee}$ is defined by $N^{\vee}_i = (N_{-i})^\ast$.  It is
known (see, for example, \cite[(4)]{co}) that an object $N$ of $\ml$ is
injective if and only if it is free, that is,
\[
N \cong \Lambda(a_1) \oplus \ldots \oplus \Lambda(a_m)
\]
or
\[
N \cong \Lambda^{\vee}(a_1) \oplus \ldots \oplus \Lambda^{\vee}(a_m)
\]
for some integers. Note that the isomorphisms 
$\Lambda^p V \to {\Hom}_k(\Lambda^{n+1-p} V, \Lambda^{n+1} V)$, 
$\omega \mapsto (-1)^p(- \wedge \omega)$, define an isomorphism
$\Lambda \cong \Lambda^{\vee} (-n-1)$ in $\ml$.

We let $\Lambda_+ \subset \Lambda$ denote the ideal generated by $V$.  Then
$\Lambda/\Lambda^{i+1}_+$ is the module $k \oplus V \oplus \ldots \oplus
\Lambda^i V$.  It is easy to show (see, for example, \cite{co},(4)(i)) that 
\vskip5mm

\subsection{Lemma}\label{determ}
\textit{  For any $N\in \Ob (\ml)$ and any integer $a$, the natural map
\[
\Hom_{\ml}\bigl(N, \Lambda^{\vee}(a)\bigr) \lra \Hom_k \bigl( N_{-a},\, 
\Lambda^{\vee}(a)_{-a}\bigr)
\]
is bijective.}
\vskip5mm

\subsection{The BGG--functor}\label{bggf}
To any $N \in \Ob(\ml)$ one associates a bounded complex $\L(N) \in
\Ob{\C}^b(\Coh \P)$ defined by $\L(N)^p = N_p \otimes \ko(p)$ 
with differential $N_p \otimes \ko(p) \to N_{p+1} \otimes \ko(p+1)$ induced 
by $N_p \otimes V \to N_{p+1}$, thus defining the BGG--functor
\[
\L : \ml \lra {\C}^b(\Coh \P)
\]
between $\ml$ and the category of bounded complexes of coherent
sheaves on $\P$.   The complex $\L(\Lambda^{\vee})$ is the standard Koszul
complex
\[
0 \to \Lambda^{n+1} V^\ast \otimes \ko(-n-1)\to \dots \to \Lambda^2 V^\ast
\otimes \ko(-2)  \to V^\ast \otimes \ko(-1) \to \ko \to 0,
\]
and the complex $\L((\Lambda/\Lambda^{i+1}_+)^{\vee})$ is the truncation
\[
0 \lra \Lambda^iV^\ast \otimes \ko (-i) \lra \dots \lra V^\ast \otimes \ko(-1)
\lra \ko \lra 0
\]
of the Koszul complex $\L(\Lambda^{\vee})$ at $-i$.

Given a complex $K^\bullet \in \Ob \C(\Lambda)$, one
obtains a double complex $X^{\bullet\bullet}$ of coherent sheaves on $\P$ with
$X^{pq} = \L(K^p)^q$ and from this the simple complex
\[
\L(K^\bullet) = \text{s}(X^{\bullet\bullet})
\]
of the double complex, which is a complex of quasi--coherent sheaves.   If
$K^\bullet$ is bounded, then $\L(K^\bullet)$ is a bounded complex of coherent
sheaves.   In this way, the functor $\L$ can be extended to functors : 
$\C(\Lambda) \to \C(\P)$, $\K(\Lambda) \to \K(\P)$, 
${\D}^+(\Lambda) \to \D(\P)$, and the last one restricts to a functor 
${\D}^b(\Lambda) \to {\D}^b(\Coh\P)$.     

For later use we need the following lemma. For a short standard proof 
see \cite{co},(5). 
\vskip5mm

\subsection{Lemma}\label{acyc1}
\textit{ Let $N^\bullet \in \Ob \C^-(\Lambda)$ be a complex of free 
(= injective) objects, bounded to the right. 
Then the complex $\Ll (N^\bullet)$ is acyclic.
If, in addition, $N^\bullet \in \Ob \C^b(\Lambda)$ is bounded, then 
$\Ll (N^\bullet)$ is bounded and consists of finite direct sums 
of line bundles.}
\vskip5mm

\subsection{Definitions for complexes} \label{def1} 

(i)  Given $N \in \Ob(\ml)$ and an integer $d$, we denote by
  $N_{\le d}$, resp.\ $N_{\ge d}$, the $\Lambda$--submodule of $N$ generated
  by its homogeneous elements of degree $\le d$, resp.\ $\ge d$.  Then
  $(N_{\ge d})_i = N_i$ for $i \ge d$ and $= 0$ for $i < d$.

\noindent If $N \xrightarrow{f} L$ is a homomorphism in $\ml$ such that
$\im(f) \subset L\cdot \Lambda_+$, then $f(N_{\le d}) \subset L_{\le d-1}$.
\vskip5mm

(ii) We let $\ki$ denote the full subcategory of $\ml$ consisting of
  injective (= free) objects.
If $I^\bullet \in \Ob \C(\ki)$, we may write each module as a finite sum
\[
I^p = \bigoplus_i H^i_{p-i} \otimes \Lambda^{\vee} (p-i),
\]
where $H^i_{p-i}$ are finite dimensional $k$--vector spaces, the coefficient
spaces of the modules $\Lambda^{\vee} (p-i)$.   We say that $I^\bullet$ is
\textbf{minimal} if $\im(d^p_I) \subset I^{p+1} \cdot \Lambda_+$ for any
$p$.   This is equivalent to
\[
d^p_I \bigl( H^i_{p-i} \otimes \Lambda^{\vee} (p-i)\bigr) \subset \bigoplus_{j
  \le i} H^j_{p+1-j} \otimes \Lambda^{\vee}(p+1-j)
\]
for any $p$, $i \in \Z$.  This shows that, if $I^\bullet$ is minimal, 
the filtration on $I^p$, given by
\[
F_j I^p := \bigoplus_{i\le j} H^i_{p-i} \otimes \Lambda^{\vee}(p-i)
\]
for any index $p$, defines a subcomplex  
$F_j I^\bullet = (F_j I^p)_{p \in \Z}$  of $I^\bullet$. Alternatively, 
$F_j I^p = I^p_{\le j-p-n-1}$. This filtration has been already 
considered by G.~Fl{\o}ystad \cite{fl1}. 
The associated graded complex $\gr_F(I^\bullet) = \underset{j}{\oplus} F_{j}
I^\bullet/F_{j-1}I^\bullet$ is the ``linear part'' of $I^\bullet$ as defined in
\cite{efs}. 
\vskip5mm

(iii) For any $K^\bullet \in \Ob \C(\Lambda)$ we define another 
subcomplex 
$\sigma K^\bullet$ of $K^\bullet$ by $\sigma K^p = K^p_{\ge - p}$.   
If $I^\bullet \in \Ob \C(\ki)$ has the terms 
written as in (ii) above, then
\[
\sigma I^p = \bigoplus_i H^i_{p-i} \otimes
(\Lambda/\Lambda_+^{i+1})^{\vee} (p-i).
\]
   It follows from the definition of the functor $\L$ that 
\[
\L(\sigma K^\bullet) = \sigma^{\ge 0} \L(K^\bullet),
\]
where $\sigma^{\ge 0}$ denotes the standard truncation of a complex killing
the terms of degree $< 0$.
\vskip5mm

\subsection{HT--complexes}\label{ht}
A minimal bounded complex $G^\bullet \in \Ob \C^b(\ki)$ will be called an
HT--complex if $F_{n-1} G^\bullet = G^\bullet$ and $F_0 G^\bullet =0$.   
The modules of such a complex can be written as
\[
G^p = \bigoplus_{0 < i < n} H^i_{p-i} \otimes \Lambda^{\vee}(p-i).
\]
\vskip5mm

\subsection{Lemma}\label{acyc2}
\textit{For any HT--complex $G^\bullet$ the complex $\Ll (G^\bullet)$ is a 
bounded acyclic complex of finite direct sums of line bundles, and
$\Ll (G^\bullet)^p=0$ for $p < -n$ or $p \ge n$.}  

The first statement follows from \ref{acyc1}, the second from the special type
of the modules $G^p$.
\vskip5mm

\subsection{Remark}\label{rht}
The sheaves of the complex $\L(G^\bullet)$ can be written as
\[
\L(G^\bullet)^s = \bigoplus_{0<i<n} \bigoplus_d H^i_d \otimes \Lambda^{i-s}
V^\ast \otimes \ko(-i-d+s)
\]
and so are precisely the terms of the acyclic complexes considered in
\cite{tr2}.   On the other hand, in \cite{ho2} G.~Horrocks employed
filtrations of bundles $E \oplus L$, $L$ a direct sum of line bundles, which
lead to similar complexes. We let
\[
\kz^0 \L(G^\bullet)
\]
denote the kernel of the differential $\L(G^\bullet)^0 \to \L(G^\bullet)^1$.
One can verify that $\kz^0 \L(G^\bullet)$ is a successive extension of the
sheaves $H^i_d \otimes \Omega^i(-d)$ and that $\H^i (\kz^0 \L(G^\bullet)(d)) 
\cong H_d^i$ for $0 < i < n$ and any $d$. For the latter see 
Corollary \ref{cohcor}.
\vskip1cm
\section{The stable category of vector bundles}\label{section2}

We denote by $\kh$ the full subcategory of the homotopy category $\K(\Lambda)$
consisting of HT--complexes.

Let, furthermore, $\kv$ denote the full subcategory of $\Coh\P(V)$ consisting 
of locally
free sheaves (vector bundles) and $\kp$ the full subcategory of $\kv$
consisting of finite direct sums of line bundles.   The \textbf{stable
  category} is the category $\kv/\kp$ which has the same objects as $\kv$ but
the sets $\Hom$ in $\kv/\kp$ are defined as 
\[
\underline{\Hom}_\ko (E,E^\prime) := \Hom_\ko(E, E^\prime)/S(E, E^\prime),
\]
where $S(E, E^\prime)$ is the subspace of homomorphisms which factorize
through an object of $\kp$.   With this definition, a morphism $E
\xrightarrow{f} E^\prime$ is an isomorphism in $\kv/\kp$ if and only if there
are objects $P$ and $P^\prime$ in $\kp$ such that $f$ factorises as
\[
E \hookrightarrow E \oplus P \overset{\approx}{\lra} E^\prime \oplus
P^\prime \twoheadrightarrow E^\prime
\]
with an isomorphism in the middle.  We need the following result of Horrocks:

\subsection{Lemma}\label{lehorr}
\textit{
  Let $E \xrightarrow{f} E^\prime$ be a homomorphism of vector bundles on
  $\P(V)$.   If the induced homomorphisms ${\Hh}^if(d)$ between ${\Hh}^iE(d)$ 
  and
  ${\Hh}^iE^\prime(d)$ are isomorphisms for all $0 < i < n$ and all $d$, then
  there exist direct sums $P$ and $P^\prime$ of line bundles such that $f$
  decomposes as
\[
E \hookrightarrow E \oplus P \overset{\approx}{\lra} E^\prime \oplus P^\prime
\twoheadrightarrow E^\prime\,.
\]  
}

This lemma is Theorem 7.5 in \cite{ho1}.   There, a short direct proof is
given, independent of the other results in \cite{ho1}.   For the conclusion of
Horrock's proof, it is convenient to consider the mapping cone of a certain
quasi--isomorphism, arising naturally from his arguments. 

\subsection{Main result.}\label{main}

The main result of this paper says that the functor
\[
G^\bullet \longmapsto \kz^0 \L(G^\bullet)
\]
induces an equivalence between the homotopy category $\kh$ of HT--complexes
and the stable category $\kv/\kp$.

This re--proves the main result of \cite{tr2}, saying that any vector bundle
$E$ on $\P(V)$ is stably equivalent to a bundle $\kz^0 \L(G^\bullet)$.

The fact that the functor $\kh \to \kv/\kp$ is fully faithful follows from
Theorem \ref{homun} below, and the fact that it is essentially surjective from
Theorem \ref{stabism}.  Moreover, we shall prove that if $I^\bullet
\xrightarrow{f^\bullet} J^\bullet$ is a morphism of minimal complexes 
of injective
objects of $\ml$, then $f^\bullet$ is already an isomorphism if 
it is a homotopy
equivalence, see Remark \ref{hisom}.  This will imply that the (usual)
isomorphism classes of HT--complexes are in bijection with the stable
isomorphism classes of vector bundles on $\P(V)$.
\vskip1cm
\section{Auxiliary equivalences}

\subsection{Lemma}\label{trequi}
\textit{Let $I^\bullet \in \Ob \C^b(\ki)$ be a bounded injective complex in
  $\ml$ and suppose that $F_0 I^\bullet = 0$.   Then, for any bounded complex
  $K^\bullet \in \Ob \C^b(\Lambda)$, the maps
\[
\Hom_{\C(\Lambda)} (K^\bullet, I^\bullet) \overset{\approx}{\lra}
\Hom_{\C(\Lambda)} (\sigma K^\bullet, I^\bullet)
\]
and
\[
\Hom_{\K(\Lambda)} (K^\bullet, I^\bullet) \overset{\approx}{\lra}
\Hom_{\K(\Lambda)} (\sigma K^\bullet, I^\bullet)
\]
are isomorphisms.}

\begin{proof}
Let $\Hom^\bullet(K^\bullet, I^\bullet)$ denote the complex defined by 
$$\Hom^q(K^\bullet,I^\bullet)=\underset{p}{\prod}\Hom_{\ml}(K^p, I^{p+q})$$
with differentials $d^q$ sending a tuple $(a^p)_p$ to  
$(a^{p+1} \circ d^p_K+(-1)^{q+1}d^{p+q}_I \circ a^p)_p$. Then 
$\Hom_{\C(\Lambda)}(K^\bullet, I^\bullet) = \Ker d^0$, and the morphisms 
$K^\bullet\to I^\bullet$ which are homotopically equivalent to zero 
are those of
$\im d^{-1}$. The statements of the lemma then follow if the maps
\[
\Hom^i(K^\bullet, I^\bullet) \overset{\approx}{\lra} \Hom^i(\sigma K^\bullet,
I^\bullet)
\]
are isomorphisms for $i = -1, 0, 1$. But this follows easily from Lemma 
\ref{determ}.
\end{proof}

\subsection{Lemma}\label{Lequi}

\textit{Let $I^\bullet \in \Ob \C^b(\ki)$ and $K^\bullet \in 
\Ob \C^b (\Lambda)$.   Then the map
\[
\Hom_{\K(\Lambda)} (K^\bullet, I^\bullet) \overset{\approx}{\lra} \Hom_{\K(\P)}
\bigl(\Ll (K^\bullet), \Ll (I^\bullet)\bigr)
\]
is an isomorphism.}

\begin{proof}
  The functor $\Hom_{\K(\Lambda)} (-, I^\bullet)$ maps quasi--isomorphisms in
  $\K^+(\Lambda)$ to isomorphisms because $I^\bullet$ consists of injective
  objects, and $\L$ maps short exact sequences to semi--split short exact
  sequences.  ($0 \to\ka^\bullet \to \kb^\bullet \to \kc^\bullet \to 0$ is
  called semi--split if each sequence $0 \to \ka^p \to \kb^p \to \kc^p \to 0$
  is split exact.)

Using the same kind of argument as at the beginning of the proof of \cite{co},
(6)(b), the proof can be reduced to the case where $K^\bullet =
\underline{k}$ and then to the case when $I^\bullet = \T^{-p}
\Lambda^{\vee}(p-i)$.   In this last case, both sides are $0$, except when $p =
i = 0$, and, if $p = i = 0$, the morphism of the lemma is clearly  
an isomorphism.
\end{proof}

The next general fact should be well--known.   We include a sketch of proof
for lack of a reference.

\subsection{Lemma}\label{extcon}
\textit{
Let $\ka$ be an abelian category with sufficiently many injective objects and
let $X^\bullet \in \Ob \C^-(\ka)$ and 
$Y^\bullet \in \Ob \C^+(\ka)$ be bounded
complexes to the right, resp.\ left.   If $\Ext^{p-q}(X^p,Y^q) = 0$ 
for all $p > q$, then the canonical map
\[
\Hom_{\K(\ka)} (X^\bullet, Y^\bullet) \lra \Hom_{\D(\ka)} (X^\bullet,
Y^\bullet)
\]
is surjective, and if $\Ext^{p-q-1}(X^p,Y^q) = 0$ for all $p > q+1$, 
then it is injective.  }

\begin{proof}
  There is a quasi--isomorphism $v^\bullet : Y^\bullet \to J^\bullet$ with
  $J^\bullet$ bounded to the left and consisting of injective objects.   It is
  well--known that then the natural homomorphism
\[
\Hom_{\K(\ka)} (X^\bullet, J^\bullet) \overset{\approx}{\lra} \Hom_{\D(\ka)}
(X^\bullet, J^\bullet)
\]
is an isomorphism for any complex $X^\bullet$ (any quasi--isomorphism 
$t^\bullet : J^\bullet \to Z^\bullet $ has a left inverse in $\K(\ka)$ 
because, $\Con(t^\bullet )$ being exact, 
$\Hom_{\K(\ka)} (\T^{-1}\Con(t^\bullet ), J^\bullet ) = 0$ ). 
Because also $\Hom_{\D(\ka)} (X^\bullet, Y^\bullet) \cong \Hom_{\D(\ka)}
(X^\bullet, J^\bullet)$, the statement of the lemma is equivalent to the
surjectivity, resp.\ injectivity, of
\begin{equation}
\Hom_{\K(\ka)} (X^\bullet, Y^\bullet) \lra \Hom_{\K(\ka)} (X^\bullet,
J^\bullet)\,. \tag{$*$} 
\end{equation}
In order to verify this, we consider the mapping cone $C^\bullet =
\Con(v^\bullet)$ of $v^\bullet$, which is exact.  We obtain the 
exact sequence
\[
\Hom_{\K(\ka)} (X^\bullet, {\T}^{-1}C^\bullet)\to 
\Hom_{\K(\ka)} (X^\bullet, Y^\bullet)
\to \Hom_{\K(\ka)} (X^\bullet, J^\bullet)\to 
\Hom_{\K(\ka)}(X^\bullet, C^\bullet)\,.
\]
The surjectivity of (*) follows from 
$\Hom_{\K(\ka)} (X^\bullet, C^\bullet) = 0$,
which we prove next.   Because $C^q = Y^{q+1} \oplus J^q$, the assumption of
the lemma for surjectivity implies that $\Ext^{p-q}(X^p, C^{q-1}) = 0$ 
for $p > q$.   Because $C^\bullet$ is exact, there are the exact sequences
\[
0 \lra Z^{q-1} (C^\bullet) \lra C^{q-1} \lra Z^q(C^\bullet) \lra 0\,.
\]
Using these, it follows that
$\Ext^1(X^p, Z^{p-1}(C^\bullet))$ injects into 
$\Ext^2(X^p, Z^{p-2}(C^\bullet))$ which injects into 
$\Ext^3(X^p, Z^{p-3}(C^\bullet))$ and so on.   But 
$Z^q(C^\bullet)\subseteq C^q$ and $C^q = 0$ for $q<<0$, hence 
$\Ext^1(X^p, Z^{p-1}(C^\bullet)) = 0$.   
This implies that any morphism $X^p \to Z^p (C^\bullet)$ can be
lifted to a morphism $X^p \to C^{p-1}$.   This can finally be used to verify
that any morphism $X^\bullet \to C^\bullet$ is homotopic to $0$, using
descending induction and the assumption that $X^\bullet$ is bounded above.   
This proves that (*) is surjective.   The condition for injectivity implies, 
in the same way, that $\Hom_{\K(\ka)} (X^\bullet, {\T}^{-1}C^\bullet ) = 0$.
\end{proof}
\vskip1cm

\section{Uniqueness}\label{section4}

We are now in position to prove that an HT--complex $G^\bullet$ is determined
by the stable isomorphism class of its associated sheaf 
$\kz^0 \L(G^\bullet)$ up to homotopy equivalence.
More precisely, there is the 

\subsection{Theorem}\label{homun}
\textit{ Let \mbox{$G^\bullet$} and \mbox{$G^{\prime\bullet}$} be
  HT--complexes and let \mbox{$\ke = \kz^0 \Ll (G^\bullet)$} and 
\mbox{$\ke^\prime  = \kz^0 \Ll (G^{\prime\bullet})$} 
be the associated locally free sheaves.
  Then the map $f^{\bullet}  \mapsto \kz^0 \Ll (f^\bullet)$ from
\[
\Hom_{\C(\Lambda)} (G^\bullet, G^{\prime\bullet}) \lra \Hom_\ko (\ke,
\ke^\prime)
\]
induces an isomorphism
\[
\Hom_{\K(\Lambda)} (G^\bullet, G^{\prime\bullet}) \overset{\approx}{\lra}
\underline{\Hom}_\ko (\ke, \ke^\prime).
\]
}
\begin{proof}
  By Lemma \ref{trequi} we have the isomorphisms
\begin{align*}
  \Hom_{\C(\Lambda)} (G^\bullet, G^{\prime\bullet}) &\cong 
\Hom_{\C(\Lambda)}(\sigma G^\bullet, G^{\prime\bullet})\\
  \Hom_{\K(\Lambda)} (G^\bullet, G^{\prime\bullet}) &\cong 
\Hom_{\K(\Lambda)}(\sigma G^\bullet, G^{\prime\bullet})
\end{align*}
and by Lemma \ref{Lequi} the isomorphism
\[
\Hom_{\K(\Lambda)} (\sigma G^\bullet, G^{\prime\bullet}) \cong \Hom_{\K(\P)}
\bigl(\L(\sigma G^\bullet), \L(G^{\prime\bullet})\bigr).
\]
Recalling $\L(\sigma G^\bullet) = \sigma^{\ge 0} \L(G^\bullet)$, 
this is a right
resolution of $\ke$ by direct sums of line bundles, of length $n-1$, similarly
$\L(G^{\prime\bullet})$.   We, thus, obtain resolutions by direct sums of 
line bundles 
\begin{align*}
  \ke^\ast & \longleftarrow \L(\sigma G^\bullet)^\ast\\
\ke^{\prime\ast} & \longleftarrow \L(\sigma G^{\prime\bullet})^\ast\,,
\end{align*}
of length $n-1$ and, therefore, induced free resolutions over the 
symmetric algebra $S$ 
\begin{align*}
  \Gamma_\ast \ke^\ast & \longleftarrow \Gamma_\ast\bigl(\L(\sigma
  G^\bullet)^\ast\bigr)\\
\Gamma_\ast \ke^{\prime\ast} & \longleftarrow \Gamma_\ast\bigl(\L(\sigma
  G^{\prime\bullet})^\ast\bigr)\,.
\end{align*}
Given a homomorphism $\ke \xrightarrow{f} \ke^\prime$ one can construct a
homomorphism
\[
\Gamma_\ast\bigl(\L(\sigma G^{\prime\bullet})^\ast\bigr) \lra \Gamma_\ast
\bigl(\L(\sigma G^\bullet)^\ast\bigr)
\]
over $\Gamma_\ast(f^{\ast})$ and, thus, a homomorphism
\[
\UseComputerModernTips
\xymatrix{
\ke \ar[r]\ar[d]_f & \L(\sigma G^\bullet)\ar[d]^{{\phi}^\bullet}\\
\ke^\prime \ar[r] & \L(\sigma G^{\prime\bullet})}
\]
of the right resolutions over $f$.   Now ${\phi}^\bullet$ 
can be extended by $0$ into
the left resolution $\sigma^{<0}\L(G^{\prime\bullet})$ of $\ke^\prime$ 
to obtain a
homomorphism $\L(\sigma G^\bullet) \to \L(G^{\prime\bullet})$, hence a 
homomorphism $f^\bullet : \sigma G^\bullet \to G^{\prime\bullet}$.     
This proves that the
map of the theorem is surjective.

If $\sigma G^\bullet \to G^{\prime\bullet}$ is homotopic to zero, then the
induced homomorphism $\ke \to \ke^\prime$ factorizes through
$\L(G^{\prime\bullet})^{-1}$ and, thus, becomes zero in the stable category.
Let, conversely, $\sigma G^\bullet \xrightarrow{f^\bullet} G^{\prime\bullet}$
define a homomorphism $\ke \xrightarrow{f} \ke^\prime$, which factorizes
through a direct sum $\kl$ of line bundles.   Because $\sigma^{<0}
\L(G^{\prime\bullet})$ is a left resolution  
of $\ke^\prime$ by direct sums of line bundles, of length $n-1$,
we obtain a free resolution over the symmetric algebra $S$ 
\[
\Gamma_\ast \sigma^{<0}\L(G^{\prime\bullet}) \lra \Gamma_\ast \ke^\prime
\]
and, then, $\Gamma_\ast \kl \to \Gamma_\ast \ke^\prime$ factorizes through
$\Gamma_\ast \L(G^{\prime\bullet})^{-1}$.    Using the dual
resolutions again, one can construct from this a homotopy 
$\Gamma_\ast (\L(f^\bullet)^\ast ) \sim 0$, hence 
$\L(f^\bullet)^{\ast}  \sim 0$, hence 
$\L(f^\bullet) \sim 0$. It follows that $f^\bullet \sim 0$. 
  This proves that
\[
\Hom_{\K(\Lambda)} (G^\bullet, G^{\prime\bullet}) \to \underline{\Hom}_\ko
(\ke, \ke^\prime)
\]
is also injective.
\end{proof}
\vskip5mm

\subsection{Remark}\label{hisom}

Let $I^\bullet$ and $J^\bullet$ be minimal complexes in  $\C(\ki)$ and let
$I^\bullet \xrightarrow{f^\bullet} J^\bullet$ be a morphism of complexes.   
If $f^\bullet$ is a homotopy equivalence 
(i.e., an isomorphism in $\K(\Lambda)$), 
then it is already an isomorphism in $\C(\Lambda)$.

\begin{proof}
Because all modules $I^p$ and $J^p$ are free, we may write again $I^p =
\underset{i}{\oplus} H^i_{p-i} \otimes \Lambda^{\vee}(p-i)$ and $J^p
=\underset{i}{\oplus} K^i_{p-i} \otimes \Lambda^{\vee}(p-i)$.   
$f^\bullet$ induces isomorphisms 
\[
\Hom_{\K(\Lambda )}(\underline{k}, {\T}^pI^\bullet (i-p)) 
\overset{\approx}{\lra} 
\Hom_{\K(\Lambda )}(\underline{k}, {\T}^pJ^\bullet (i-p)) 
\]
$\forall p,\  i \in \Z$, where $\underline{k}=\Lambda/\Lambda_+$.   
Since $I^\bullet$ is minimal, 
\[
\Hom_{\K(\Lambda )}(\underline{k}, {\T}^pI^\bullet (i-p)) \cong 
H^i_{p-i}\otimes (\Lambda^{\vee} )_0
\]
and, analogously, for the other term. One deduces that the component of 
$f^p_{i-p}$ mapping 
$H^i_{p-i}\otimes \Lambda^{\vee} (p-i)_{i-p}$ to 
$K^i_{p-i}\otimes \Lambda^{\vee} (p-i)_{i-p}$ is an isomorphism, hence the 
component of $f^p$ mapping $H^i_{p-i}\otimes \Lambda^{\vee} (p-i)$ to 
$K^i_{p-i}\otimes \Lambda^{\vee} (p-i)$ is an isomorphism 
$\forall p,\  i \in \Z$. 
It follows that each $f^p$ is an isomorphism (because the matrix defining 
it is triangular).
\end{proof}

The following isomorphism will be needed later. Its corollary is the 
result 3.7 in \cite{tr2}.
\vskip5mm

\subsection{Lemma}\label{cohom}
\textit{Let $I^\bullet \in \Ob{\C}^b(\ki)$ and let $\underline{k} =
  \Lambda/\Lambda_+$.   Then the natural homomorphism
\[
\Hom_{\K(\Lambda)} \bigl(\underline{k},(\sigma I^\bullet)(m)\bigr)
\overset{\approx}{\lra} \Hom_{\D(\P)} \bigl(\Ll (\underline{k}), 
\Ll \bigl((\sigma I^\bullet)(m)\bigr)\bigr)
\]
is an isomorphism for $0 < m < n$. }

\begin{proof}
$(\sigma I^p) (m)_0 = I^p(m)_0$ for $p\ge -1$ and $m > 0$ and so
\[
\Hom_{\K(\Lambda)} \bigl(\underline{k}, (\sigma I^\bullet)(m)\bigr) \cong
\Hom_{\K(\Lambda)}\bigl(\underline{k}, I^\bullet(m)\bigr)\,.
\]
As in \ref{def1}, we have $\L\bigl((\sigma I^\bullet)(m)\bigr) = \sigma^{\ge 
-m}\L\bigl(I^\bullet(m)\bigr)$.   Therefore, also
\[
\Hom_{\K(\P)} \bigl(\L(\underline{k}), 
\L\bigl((\sigma I^\bullet)(m)\bigr)\bigr) \cong
\Hom_{\K(\P)} \bigl(\L(\underline{k}), \L\bigl(I^\bullet(m)\bigr)\bigr)
\]
for $m > 0$.   Now Lemma \ref{Lequi} applied to the free (= injective) complex
$I^\bullet(m)$ implies that the natural homomorphism
\[
\Hom_{\K(\Lambda)} \bigl(\underline{k}, (\sigma I^\bullet)(m)\bigr)
\overset{\approx}{\lra} 
\Hom_{\K(\P)}\bigl(\L(\underline{k}), 
\L\bigl((\sigma I^\bullet)(m)\bigr)\bigr)
\]
is an isomorphism.   Because $\L(\underline{k}) = \ko $ 
, the conditions of Lemma \ref{extcon} are satisfied for $m < n$ 
so that the last $\Hom_{\K(\P)}$ is isomorphic to $\Hom_{\D(\P)}$.
\end{proof}
\vskip5mm

\subsection{Corollary}\label{cohcor}
\textit{If $G^\bullet $ is an HT--complex then 
${\Hh}^m(\kz^0 \Ll (G^\bullet)(p-m)) \cong H^m_{p-m}$ for $0<m<n$ 
and all $p\in\Z.$}

\begin{proof}
It follows from the definition of the functor $\L$ that $\L\bigl(\T^p(\sigma
G^\bullet)(m-p)\bigr) = \T^m \L(\sigma G^\bullet)(p-m)$.   
Applying the previous
lemma to $I^\bullet = {\T}^pG^\bullet (-p)$ we obtain 
(because $\sigma I^\bullet = {\T}^p(\sigma G^\bullet)(-p)$) an isomorphism
\[
\Hom_{\K(\Lambda)} \bigl(\underline{k}, \T^p(\sigma G^\bullet)(m-p)\bigr) 
\cong
\Hom_{\D(\P)} \bigl(\L(\underline{k}), \T^m \L(\sigma G^\bullet)(p-m)\bigr)\,.
\]
Now it is easy to verify that the left--hand side is isomorphic to $H^m_{p-m}$.
Because $\L(\sigma G^\bullet)$ is quasi--isomorphic to the trivial complex $\ke
= \kz^0 \L(G^\bullet)$, the right--hand side is isomorphic to $\Hom_{\D(\P)}
\bigl( \ko, \T^m \ke(p-m)\bigr) = \Ext^m \bigl(\ko, \ke(p-m)\bigr) 
\cong \H^m\ke(p-m)$. 
\end{proof}
\vskip5mm

\subsection{Remark}\label{multip}
One can show, moreover, that the linear part 
$H^m_{p-m} \otimes {\Lambda}^{\vee} (p-m) \to 
H^m_{p+1-m} \otimes {\Lambda}^{\vee} (p+1-m)$ 
of the differential $d^m_G$ is induced (up to sign) by the multiplication 
map $\H^m{\ke}(p-m) \otimes V^\ast \to \H^m{\ke}(p+1-m)$. 

Indeed, let $I^\bullet \in \Ob{\C}^b(\ki)$ and $0 < m < n$. We have seen, 
in the proof of \ref{cohom}, that one has a canonical isomorphism 
\[
\Hom_{\K(\Lambda )}\bigl(\underline{k}, I^\bullet (m)\bigr) \cong 
\Hom_{\D(\P)}\bigl(\L(\underline{k}), 
\L\bigl( (\sigma I^\bullet )(m)\bigr)\bigr)
\]
and, in the same way, one proves that there is a canonical isomorphism 
\[
\Hom_{\K(\Lambda )}\bigl({\T}^{-1}\underline{k}(1), I^\bullet (m)\bigr) 
\cong 
\Hom_{\D(\P)}\bigl(\L\bigl({\T}^{-1}\underline{k}(1)\bigr), 
\L\bigl( (\sigma I^\bullet )(m)\bigr)\bigr)\,.
\]
Now, one may use the arguments from the proof of \cite{co},(10)(b). 
\vskip5mm

For the proof of Theorem \ref{stabism} below we need also the following 
lemma which is \cite{bgg}, Remark 3 after Theorem 2, in the form stated 
and proved in \cite{co},(6)(b). 
\vskip5mm

\subsection{Lemma}\label{cohom2}
\textit{Let $I^\bullet \in \Ob \C(\ki)$ be an acyclic complex.  
Then the canonical morphism 
\[
\Hom_{\K(\Lambda)}(\underline{k},I^\bullet)
\overset{\approx}{\lra} \Hom_{\D(\P)}(\Ll (\underline{k}), \Ll (I^\bullet))
\]
is an isomorphism.}
\vskip1cm

\section{Comparison with Tate resolutions}\label{section5}

It was shown in \cite{tr2} that any stable isomorphism class of a 
vector bundle $\ke$ on
$\P V$ is the class of $\kz^0 \L(G^\bullet)$ of an HT--complex $G^\bullet$ with
$H^i_{p-i} = {\H}^i \ke (p-i)$ for $0 < i < n$, and that this complex is unique
up to isomorphisms.   On the other hand, the main result of \cite{bgg} was
that there is a module $N \in \Ob(\ml)$, annihilated by
soc$(\Lambda) = \Lambda_{n+1}$, such that $\L(N)$ is a monad for $\ke$, that 
is, ${\kh}^0 \L(N) \cong \ke$ and ${\kh}^i \L(N) = 0$ for $i \not= 0$.  
Also in this
situation $N$ is unique up to isomorphism.   In \cite{ge}, and
recently in \cite{efs}, so called Tate resolutions of $\Lambda$--modules
associated to bundles or coherent sheaves on $\P V$ have been used in order to
improve the understanding of the result of \cite{bgg}.  A Tate resolution in
our setting is an acyclic minimal complex $I^\bullet$ in $\C(\ki)$ such that 
$N \cong Z^0(I^\bullet) = \Ker(I^0 \to I^1)$.   We say that $I^\bullet$ is a
Tate resolution of $\ke$.   This agrees with the definition in 
\cite{efs}, Section 4,  as one may verify by using \cite{co},(10). The main
result in \cite{efs} is that
\[
I^p \cong \underset{0 \le i \le n}{\oplus} \H^i \ke (p-i) \otimes \Lambda^{\vee}
(p-i)
\]
and that the linear part
\[
\H^i \ke(p-i) \otimes \Lambda^{\vee} (p-i) \lra \H^i \ke(p+1-i) \otimes
\Lambda^{\vee} (p+1-i)
\]
of the differential $I^p \to I^{p+1}$ is induced by the multiplication map 
(also
valid for hyper-cohomology of a complex $\kf^\bullet$ of coherent sheaves).
This was stated in \cite{bgg}, Remark 3 after Theorem 2, 
and completely proved
in \cite{efs}, see also \cite{co} for an exposition of these results in a
spirit closer to that of \cite{bgg}.

The relation between HT--complexes $G^\bullet$ and the Tate resolutions for a
bundle $\ke$ is described in the following theorem and its corollary.
\vskip5mm

\subsection{Theorem}\label{stabism} 
\textit{ Let $\ke$ be a vector bundle on $\P(V)$
and $I^\bullet$ a Tate resolution of $\ke$.   Then $\ke$ is stably isomorphic
to $\kz^0 \Ll (F_{n-1} I^\bullet / F_0 I^\bullet)$.}

\begin{proof}
  (1) We may assume that $I^p = \underset{i}{\oplus} H^i_{p-i} \otimes
      \Lambda^{\vee} (p-i)$ with $H^i_d = \H^i \ke(d)$.   It follows that
\[
G^\bullet := F_{n-1} I^\bullet / F_0 I^\bullet
\]
is an HT--complex.   $G^\bullet$ is a subcomplex of $I^\bullet / F_0
I^\bullet$ but not necessarily of $I^\bullet$.   However, $\sigma G^\bullet$
is a subcomplex of $I^\bullet$ because 
\[
\sigma G^p = \underset{0<i<n}{\bigoplus} H^i_{p-i} \otimes 
(\Lambda /{\Lambda}^{i+1}_+)^{\vee} (p-i)
\]
and any morphism 
$(\Lambda /{\Lambda}^{i+1}_+)^{\vee} (p-i) \to {\Lambda}^{\vee} (p+1)$ 
is 0 by \ref{determ}. 

(2) The complex $\L(I^\bullet)$ consists of quasi--coherent sheaves which are
    at most countable direct sums of line bundles, 
    ${\kh}^0 \L(I^\bullet) \cong \ke$ and ${\kh}^i \L(I^\bullet) = 0$ 
    for $i \not= 0$ 
    by \cite{co},(6)(a). 

(3) Now let $\kf = \kz^0 \L(G^\bullet) = \kz^0 \L(\sigma G^\bullet)$.  Because
    $\L(\sigma G^\bullet)$ is a right resolution of $\kf$, the inclusion
    $\sigma G^\bullet \hookrightarrow I^\bullet$ induces a homomorphism 
    $\kf \xrightarrow{f} \ke$ via $\L(\sigma G^\bullet) \to \L(I^\bullet)$.  
    We are going to show that $f$ is the stable isomorphism of the theorem.  
    For that, we consider the diagram
\[
\UseComputerModernTips
\xymatrix{
\Hom_{\K(\Lambda)} \bigl(\underline{k}, \T^p(\sigma G^\bullet) (m-p)\bigr)
\ar[r]^a \ar[d]_b & \Hom_{\K(\Lambda)} \bigl(\underline{k}, \T^p (I^\bullet)
(m-p)\bigr)\ar[d]_c\\
\Hom_{\D(\P)} \bigl(\L(\underline{k}), \L(\T^p(\sigma G^\bullet)(m-p))\bigr)
\ar[r]^d & \Hom_{\D(\P)} \bigl(\L(\underline{k}), 
\L(\T^p(I^\bullet)(m-p))\bigr)}
\]
for any $p \in \Z$ and $0 < m< n$.   The map $a$ is an isomorphism because a
map $\underline{k} \to {\T}^p(\sigma G^\bullet)(m-p)$ can have 
its target only in
$H^m_{p-m} \otimes \underline{k}$ and, similarly, for $I^\bullet$, and because 
$I^\bullet $ is minimal.   The map $b$ is an
isomorphism by Lemma \ref{cohom}, and $c$ is an isomorphism by 
Lemma \ref{cohom2}. It follows that $d$ is an isomorphism.
But $d$ can be identified with 
${\H}^mf(p-m) : {\H}^m\kf (p-m) \to {\H}^m\ke (p-m)$.   Now Lemma \ref{lehorr} 
implies that $f$ is an isomorphism of the stable classes of $\kf$
and $\ke$. 
\end{proof}
\vskip5mm

\subsection{Corollary} \label{htisom} \textit{ Let $G^\bullet$ be any
  HT--complex, 
$\ke = \kz^0 \Ll (G^\bullet)$ its associated bundle and let
$I^\bullet$ be the Tate resolution of $\ke$.   Then
\[
G^\bullet \cong F_{n-1} I^\bullet / F_0 I^\bullet
\]
in the category $\C(\Lambda)$ of complexes.}

\begin{proof}
  Let $\kf = \kz^0 \L(F_{n-1} I^\bullet/F_0I^\bullet)$ as in the previous
  proof.   By the theorem, $\ke$ and $\kf$ are stably equivalent.   Then, by
  Theorem \ref{homun}, there is a homotopy equivalence $G^\bullet
  \xrightarrow{f^\bullet} F_{n-1} I^\bullet / F_0 I^\bullet$.   By Remark
  \ref{hisom}, $f^\bullet$ must be an isomorphism of complexes.
\end{proof}
\vskip5mm

\subsection{Remark}\label{stabism2}
Let $I^\bullet \in \Ob {\C}^b(\ki)$ be a minimal bounded complex of 
injective objects of $\ml$, and let $\ke = \kz^0\L(I^\bullet )$. Then 
$\ke$ is stably isomorphic to $\kz^0\L(F_{n-1}I^\bullet /F_0I^\bullet )$. 

Indeed, one can repeat the arguments from the proof of \ref{stabism}, 
using \ref{cohom} instead of \ref{cohom2}.
\vskip1cm

\section{HT--resolutions}

Each stable isomorphism class of a vector bundle contains a unique element,
which has no direct summand of rank 1.   The next result, which was proved in
\cite{tr2}, Section 8, shows how one can compute the invariants of this bundle
in terms of the corresponding HT--complex.

\subsection{Proposition}\label{rkdec}
\textit{
Let $G^\bullet$ be an HT--complex, with $G^p = \underset{0<i<n}{\oplus}
  H^i_{p-i} \otimes \Lambda^{\vee} (p-i)$.   Then}

  \begin{enumerate}
  \item[(a)] rank $\kz^0 \Ll (G^\bullet) = \underset{0<i<n}{\sum}
\left(\begin{smallmatrix}n\\i\end{smallmatrix}\right)  
\underset{p}{\sum} \dim(H^i_{p-i})$.

\item[(b)] \textit{If $r(p)$ denotes the rank of $G^{p-1}_{-p} \to G^p_{-p}$,
    then $\kz^0 \Ll (G^\bullet) \cong \ke \oplus \underset{p}{\oplus}
    \ko(-p)^{r(p)}$, where $\ke$ contains no direct summand of rank 1.}
  \end{enumerate}

  \begin{proof}
    (a) rank $\kz^0 \L(G^\bullet)$ is an additive function with respect to 
    $G^\bullet$, hence we may assume that 
    $G^\bullet = {\T}^{-p}{\Lambda}^{\vee} (p-i)$. In this case, 
    $\kz^0 \L(G^\bullet) = {\Omega}^i_{\P}(i-p)$. 
    
    (b) Note that the map $d^{p-1}_{G,-p}$ between the vector spaces
    $G^{p-1}_{-p}$ and $G^p_{-p}$ decomposes into
\[
\underset{0<i<n}{\oplus} H^i_{p-1-i} \otimes \Lambda^{i+1} V^\ast \lra
\underset{0<i<n}{\oplus} H^i_{p-i} \otimes \Lambda^i V^\ast\,.
\]
We have
\[
\L(G^\bullet)^{-1} = \underset{p}{\oplus}\, G^{p-1}_{-p} \otimes \ko (-p)\quad 
\text{
  and }\quad \L(G^\bullet)^0 = \underset{p}{\oplus}\, G^p_{-p} \otimes \ko(-p)\,.
\]
We are going to cancel the constant part of the differential
$d^{-1}_{\L(G)}:\L(G^\bullet)^{-1} \to \L(G^\bullet)^0$.   
Let $Z^{p-1}$ (resp., $B^p$) be the kernel (resp., image) of 
$d^{p-1}_G:G^{p-1} \to G^p$.   
In degree $-p$, we choose decompositions 
\[
Y^{p-1}_{-p} \oplus Z^{p-1}_{-p} = G^{p-1}_{-p}\quad \text{ and }   
\quad B^p_{-p} \oplus C^p_{-p} = G^p_{-p}
\]
such that $Y^{p-1}_{-p} \cong B^p_{-p}$ under $d^{p-1}_G$. 
Then the differential $\L(G^\bullet )^{-1} \to \L(G^\bullet )^0$ 
decomposes into 
\[
\left(\begin{smallmatrix}\alpha &  \beta\\ \gamma &
        \delta\end{smallmatrix}\right):\, 
\ky \oplus \kg^{-1} \longrightarrow \kb \oplus \kg^0 
\]
where the sheaves in this decomposition are defined according to the above 
decomposition, i.e., 
$\kg^{-1} = \underset{p}{\oplus}Z^{p-1}_{-p} \otimes \ko(-p)$, 
$\kg^0 = \underset{p}{\oplus}C^p_{-p} \otimes \ko(-p)$ and, similarly, 
$\ky$ and $\kb$.  
Since $\alpha$ is an isomorphism, we may consider the automorphisms  
\[
\Phi = \left(\begin{smallmatrix}\id & 0\\ \gamma^{\prime} &
    \id\end{smallmatrix}\right) 
\,:\, \L(G^\bullet )^0 \overset{\approx}{\lra} \L(G^\bullet )^0 
\quad \text{and}\quad    
\Psi = \left(\begin{smallmatrix}\id & \beta^{\prime}\\ 0 & \id\end{smallmatrix}\right)\,:\, \L(G^\bullet )^{-1} \overset{\approx}{\lra} 
\L(G^\bullet )^{-1} 
\]
where $\gamma^{\prime} = -\gamma \alpha^{-1}$ and 
$\beta^{\prime} = -\alpha^{-1}\beta $.    Then we have 
\[
\Phi \circ d^{-1}_{\L(G)} \circ \Psi =   
\left(\begin{smallmatrix}\alpha & 0\\ 0 &
    \delta^{\prime}\end{smallmatrix}\right)\,, 
\]
where ${\delta}^{\prime} = \delta - \gamma \alpha^{-1} \beta 
:{\kg}^{-1} \to {\kg}^0$. 
Since $d^{-1}_{\L(G)}$ maps $Z^{p-1}_{-p} \otimes \ko(-p)$ to 
$\oplus_{q < p} G^q_{-q} \otimes \ko(-q)$,    
${\delta}^{\prime}$ maps each $Z^{p-1}_{-p} \otimes \ko(-p)$ into
$\oplus_{q < p} C^q_{-q} \otimes \ko(-q)$.   This means that
${\delta}^{\prime}$ has no constant component. 
We have an induced acyclic complex
\[
\cdots \lra \L(G^\bullet)^{-2} \lra \kg^{-1} 
\overset{{\delta}^{\prime}}{\lra} \kg^0
  \lra \L(G^\bullet)^1 \lra \cdots 
\]
where $\L(G^\bullet )^{-2} \to \kg^{-1}$ and 
$\kg^0 \to \L(G^\bullet )^1$ are just the components of 
$d^{-2}_{\L(G)}$ and $d^0_{\L(G)}$, respectively. 
Let $\ke$ denote the image of ${\delta}^{\prime}$.   Then
\[
\ke \oplus \bigoplus_p  B^p_{-p} \otimes \ko(-p) = 
\ke \oplus \kb = 
\im (\Phi \circ d^{-1}_{\L(G)} \circ \Psi) \cong   
\im d^{-1}_{\L(G)} = \kz^0 \L(G^\bullet)\,.
\]
It remains to show that $\ke$ has no direct summand of rank 1.   Suppose there
is such a summand.   Then we have a decomposition $\ko \xrightarrow{\rho} \ke
(a) \xrightarrow{\pi} \ko$ with $\pi \circ \rho = \id$ for some integer $a$.
Because the left resolution of $\ke$ has length $n-1$, $\Gamma \kg^{-1}(a) \to
\Gamma \ke(a)$ is surjective and $\rho$ factorizes through $\kg^{-1}(a)$.   By
the same argument applied to $\ke^\ast$, also $\ke(a) \to \ko$ extends to
$\kg^0(a) \to \ko$. But this would imply that $\kg^{-1} \to \kg^0$ has a
non--zero component $\ko(-a) \to \ko(-a)$, contradicting the property of
${\delta}^{\prime}$ stated above.
  \end{proof}
\vskip5mm

\subsection{Remark}\label{minim}
By pursuing the process of cancelation of constant parts of the 
differentials, we eventually obtain from the complex $\L(G^\bullet )$ an 
acyclic complex 
\[
0 \lra \kl^{-n} \lra \cdots \lra \kl^{-1} \lra \kl^0 \lra \cdots 
\lra \kl^{n-1} \lra 0
\]
with terms 
\[
\kl^m \cong \bigoplus_p H^p(G^\bullet )_{m-p} \otimes \ko(m-p),
\]
whose differentials have 
no constant part, and such that $\ke \cong \kz^0\kl^{\bullet}$. 
One has $H^i({\Gamma}_*\kl^{\bullet}) \cong {\H}^i_*\ke$ as graded 
$S$--modules for $0 < i < n$, and $H^i({\Gamma}_*\kl^{\bullet}) = 0$ for 
$i \le 0$. 

Moreover, the ``linear part'' of $\kl^{\bullet}$ is isomorphic to 
$\L\bigl(H^\bullet (G^\bullet )\bigr)$, where $H^\bullet (G^\bullet )$ is 
the complex whose $p$th term is $H^p(G^\bullet )$ and with all the 
differentials equal to 0.  
\vskip5mm

\subsection{Example}(Eilenberg--Maclane bundles)\label{eilmac}
Let us consider the particular case when $G^\bullet $ consists of 
only one linear strand. In order to do that, it is convenient to recall 
another functor considered in \cite{bgg}. Namely, let $S$--mod be the 
category of finitely generated graded $S$--modules. One defines a functor 
R : $S$--mod $\to \C^+(\ml)$ by associating to each graded 
$S$--module $M$ a complex $\text{R}(M)$ whose $p$th term is 
$\text{R}(M)^p := M_p \otimes {\Lambda}^{\vee} (p)$ and whose $p$th 
differential $d^p_{\text{R}(M)}$ is induced by the multiplication map 
$M_p \otimes V^\ast \to M_{p+1}$. If $M$ has finite length and $0 < j < n$ 
the $G^\bullet := \text{R}(M)(-j)$ is an HT--complex. In this case, let 
$\ke$ be the vector bundle corresponding to $G^\bullet $ as in 
\ref{rkdec}, and let $\kl^{\bullet}$ be the complex obtained 
from $\L(G^\bullet )$ as in \ref{minim}. 
Using \ref{rht}, $\L(G^\bullet )^m = 0$ for $m > j$ or for 
$m < j-n-1$, hence the same is true for $\kl^{\bullet}$. It 
follows that, if 
\[
0 \to F^{-n-1} \to \dots \to F^{-1} \to F^0 \to M \to 0
\]
is a minimal free resolution of $M$ in $S$--mod, then 
${\Gamma}_*\kl^{\bullet} \cong {\T}^{-j}F^\bullet $ and $\ke$ is 
a sheafification of $\Ker(F^{-j} \to F^{-j+1})$. 
The vector bundles of this kind are called Eilenberg--Maclane bundles 
by Horrocks \cite{ho2}. 
\vskip5mm

\subsection{Remark}\label{beilin} 
Under some assumptions, the vector bundle $\ke$ from the statement of 
\ref{rkdec} is the cohomology of a Beilinson monad obtained directly 
from the HT--complex $G^\bullet $. More precisely, assume that 
$H^p(G^\bullet )_{-1-p} = 0$ for all $p < 0$, and $H^p(G^\bullet )_{-p} = 0$ 
for all $p > 0$. Let $C^\bullet $ be the complex defined by 
$C^p := \kz^0\L(G^p) \cong \underset{0 < i < n}{\oplus}H^i_{p-i} \otimes 
{\Omega}^{i-p}(i-p)$ and with the differential $d^p_C := \kz^0\L(d^p_G)$ 
for all $p \in \Z$. Then $\kh^0(C^\bullet ) \cong \ke$ and 
$\kh^i(C^\bullet ) = 0$ for $i \not= 0$. 
\vskip3mm

Indeed, we shall use the method of Eisenbud et al. \cite{efs}, (6.1) 
for deriving Beilinson monads from Tate resolutions. Let $I^\bullet $ be 
a Tate resolution of $\ke$ and consider the complex ${\tilde C}^\bullet $ 
defined by ${\tilde C}^p := \kz^0\L(I^p) \cong 
\underset{0 \le i \le n}{\oplus}\H^i\ke(p-i) \otimes {\Omega}^{i-p}(i-p)$ 
and with the differential $d^p_{\tilde C} := \kz^0\L(d^p_I)$ for all 
$p \in \Z$. Then, by \cite{efs}, (6.1)(see, also, the proof of 
\cite{co}, (12)), $\kh^0({\tilde C}^\bullet ) \cong \ke$ and 
$\kh^i({\tilde C}^\bullet ) = 0$ for $i \not= 0$. Now, consider the 
morphisms of complexes $G^\bullet \to I^\bullet /F_0I^\bullet \gets 
I^\bullet $.    
Using the notations from \ref{minim}, our additional assumptions imply that  
\[
\kl^{-1} \cong \bigoplus_{p\ge 0} H^p(G^\bullet )_{-1-p} \otimes 
\ko(-1-p)\quad \text{ and }\quad 
\kl^0 \cong \bigoplus_{p\le 0} H^p(G^\bullet )_{-p} \otimes \ko(-p),
\]
hence $\H^0\ke = 0$ and $\H^0\ke^{\ast}(-1) = 0$. It follows that 
$\H^0\ke(i) = 0$ for all $i \le 0$ and, by Serre duality, 
$\H^n\ke(i) = 0$ for all $i \ge -n$. But this implies that one gets 
isomorphisms $C^p = \kz^0\L(G^p) \overset{\approx}{\lra} 
\kz^0\L(I^p/F_0I^p) \overset{\approx}{\lla} \kz^0\L(I^p) = {\tilde C}^p$, 
for all $p \in \Z$. 
\vskip5mm

\subsection{Remark}\label{indec}
It follows from the main results of this paper that the vector bundle 
$\ke $ from the statement of \ref{rkdec} is indecomposable if and 
only if the corresponding HT--complex $G^\bullet $ is indecomposable. 
This has an amusing consequence : let $\ke$ be an indecomposable 
vector bundle on $\P_n$ and $m \in \Z$. If $\H^i\ke(m-i) = 0$ for 
$0 < i < n$ then $\H^i\ke(m^{\prime}-i) = 0$ for $0 < i < n$ and for 
either every $m^{\prime}\ge m$ or every $m^{\prime}\le m$. 

Indeed, let $G^\bullet $ be the HT--complex corresponding to the stable 
isomorphism class of $\ke$. Our assumption implies that $G^m = 0$, 
hence $G^\bullet = {\sigma}^{>m}G^\bullet \oplus 
{\sigma}^{<m}G^\bullet $. One deduces that either 
${\sigma}^{>m}G^\bullet = 0$ or ${\sigma}^{<m}G^\bullet = 0$. 
\vskip5mm

\subsection{Remark}\label{constr}
It seems unpractical to use the results from this section to construct 
vector bundles. However, one can, at least, recuperate some of the 
well--known vector bundles of small rank on projective spaces. The 
following method is essentially due to Horrocks.   Consider an HT--complex 
$G^\bullet $ with only two non--zero terms, $G^0$ and $G^{-1}$. It has 
only one non--zero differential, $d^{-1} : G^{-1} \to G^0$. Let $\ke$ be 
the vector bundle associated to $G^\bullet $ as in \ref{rkdec}. Assume 
that the conditions from \ref{beilin} are satisfied, i.e., that 
$d^{-1}_0 : G^{-1}_0 \to G^0_0$ is injective. In this case, if 
$\kl^\bullet $ is the complex from \ref{minim}, then 
\[
\kl^{-1} \cong H^0(G^\bullet )_{-1} \otimes \ko(-1)\quad \text{and}\quad 
\kl^0 \cong H^0(G^\bullet )_0 \otimes \ko \oplus 
H^{-1}(G^\bullet )_1 \otimes \ko(1)
\]
and the linear part of the differential $\kl^{-1} \to \kl^0$ is defined 
by the multiplication map $H^0(G^\bullet )_{-1} \otimes V \to 
H^0(G^\bullet )_0$. 

Let $H\subset H^0(G^\bullet )_{-1}$ (resp., $H^0(G^\bullet )_0 
\twoheadrightarrow K$) be a sub-- (resp., quotient) vector space, and 
consider the induced morphisms $\mu : H \otimes \ko(-1) \to \ke$ and 
$\varepsilon : \ke \to K \otimes \ko$. Then $\varepsilon \circ \mu = 0$ 
if and only if the composite map $H \otimes V \to H^0(G^\bullet )_0 
\twoheadrightarrow K$ is 0. Moreover, the Lemma \ref{monoepi} below gives 
conditions under which $\mu $ is a monomorphism of vector bundles and 
$\varepsilon $ is an epimorphism. If these conditions are satisfied, 
then we get a monad 
\[
0 \lra H \otimes \ko(-1) \lra \ke \lra K \otimes \ko \lra 0 
\]
whose cohomology is a vector bundle of rank $\rank(\ke) - \dim(H) - \dim(K)$.
\vskip5mm

\subsection{Lemma}\label{monoepi}
\textit{Under the assumptions of }\ref{constr}
   \begin{enumerate}
   \item[(a)] $\mu $\textit{ is a monomorphism of vector bundles if and 
only if }$H$\textit{ intersects the image of the bilinear multiplication 
map }$H^0(G^\bullet )_{-2} \times V \to H^0(G^\bullet )_{-1}$ 
\textit{ only in }0. 
\vskip2mm

   \item[(b)] $\varepsilon $\textit{ is an epimorphism if and only if 
the composite map }$H^0(G^\bullet )_{-1} \otimes V \to 
H^0(G^\bullet )_0 \twoheadrightarrow K$\textit{ induces, for every }
$0 \not= v \in V$,\textit{ a surjection } 
$H^0(G^\bullet )_{-1} \otimes v \twoheadrightarrow K$. 
\vskip2mm

   \item[(c)] \textit{Assume, furthermore, that }$d^{-1}_0$\textit{ is 
an isomorphism, hence }$H^0(G^\bullet )_0 = 0$.\textit{ Consider a 
quotient vector space }$H^{-1}(G^\bullet )_1 \twoheadrightarrow Q$
\textit{ and the induced morphism }$\pi : \ke \to Q \otimes \ko(1)$. 
\textit{ Then }$\pi$\textit{ is an epimorphism if and only if the 
subspace }$Q^\ast \subset \left(H^{-1}(G^\bullet )_1\right)^\ast $ 
\textit{ intersects the image of the bilinear multiplication map }
$\left(H^{-1}(G^\bullet )_2\right)^\ast \times V \to 
\left(H^{-1}(G^\bullet )_1\right)^\ast $\textit{ only in }0. 
\textit{ We are unable to state, in this case, a condition equivalent 
to }$\pi \circ \mu = 0$.
   \end{enumerate} 

   \begin{proof}
For (a) one uses the exact sequence 
\[
H^0(G^\bullet )_{-2} \otimes \ko(-2) \lra 
H^0(G^\bullet )_{-1} \otimes \ko(-1) \lra \ke \lra 0\,,
\]
for (b) the fact that $H^0(G^\bullet )_{-1} \otimes \ko(-1) \to \ke$ 
is an epimorphism, and for (c) the exact sequence 
\[
0 \lra \ke \lra H^{-1}(G^\bullet )_1 \otimes \ko(1) \lra 
H^{-1}(G^\bullet )_2 \otimes \ko(2) \,.
\]
    \end{proof} 

\subsection{Example}(Trautmann \cite{tr1}, Vetter \cite{ve}, 
Tango \cite{ta})\label{tango} 
Take $G^\bullet = {\Lambda}^{\vee} (-j)$ for some 
$0 < j < n$, hence $\ke = {\Omega}^j_{\P}(j)$. Consider a vector subspace  
$H\subseteq G^0_{-1} = {\Lambda}^{j+1}V^\ast $. Lemma \ref{monoepi},(a) 
tells us, in this case, that the evaluation morphism $H \otimes \ko(-1) 
\to {\Omega}^j_{\P}(j)$ is a monomorphism of vector bundles if and only if 
the $H$ intersects the image of the bilinear contraction map 
$\Lambda^{j+2}V^\ast \times V \to \Lambda^{j+1}V^\ast $ only in 0.  
It is convenient to use the identifications 
${\Lambda}^{j+1}V^\ast \cong {\Lambda}^{n-j}V$ and 
${\Lambda}^{j+2}V^\ast \cong {\Lambda}^{n-j-1}V$. 
Then the above contraction 
map can be identified (up to sign) with the exterior multiplication 
${\Lambda}^{n-j-1}V \times V \to {\Lambda}^{n-j}V$.  
 
Assume, now, that $j = n-2$. The image of the bilinear map 
$V \times V \to {\Lambda}^2V$ is the affine cone over the Pl\" ucker 
embedding of $\text{Grass}_2(V)$, hence its dimension is $2n-1$. It 
follows that the largest possible dimension of $H$ for which 
$H \otimes \ko(-1) \to {\Omega}^{n-2}_{\P}(n-2)$ is a monomorphism of 
vector bundles is 
$\left( 
  \begin{smallmatrix}
    n+1\\2
  \end{smallmatrix}\right)-(2n-1) = 
 \left(
  \begin{smallmatrix}
    n-1\\2
  \end{smallmatrix}\right)$. In this case, the rank of the 
cokernel $\kf$ is $n-1$. 

One can construct concrete examples of such 
subspaces $H$ of ${\Lambda}^2V$. For instance, if $e_0, \dots, e_n$ form a
basis of $V$, one may take for $H$ the subspace spanned by the elements
\[
w_{ij} = 
\begin{cases}
  e_i \wedge e_j - e_0 \wedge e_{i+j} & \text{for $i + j \le n$}\\
e_i \wedge e_j - e_{i+j-n} \wedge e_n & \text{for $i+j > n$}
\end{cases}
\]
where $0 < i < j < n$. Since ${\Omega}^{n-2}_{\P}(n-2)$ has a resolution 
\[
0 \to {\Lambda}^{n+1}V^\ast \otimes \ko(-3) \to 
{\Lambda}^nV^\ast \otimes \ko(-2) \to {\Lambda}^{n-1}V^\ast \otimes 
\ko(-1) \to {\Omega}^{n-2}(n-2) \to 0
\]
or, equivalently, a resolution 
\[
0 \to \ko(-3) \to V \otimes \ko(-2) \to {\Lambda}^2V \otimes \ko(-1) \to 
{\Omega}^{n-2}(n-2) \to 0
\]
it follows that $\kf$ has a resolution 
\[
0 \to \ko(-3) \to V \otimes \ko(-2) \to ({\Lambda}^2V/H) \otimes \ko(-1) 
\to \kf \to 0\,.
\]
Now, let $X_0, \dots, X_n$ be the dual basis of $V^\ast $. Choosing a 
convenient (and rather obvious) basis for the subspace 
$({\Lambda}^2V/H)^\ast $ of ${\Lambda}^2V^\ast $ (for the particular $H$ 
considered above), one sees easily that $\kf^\ast $ is the 
vector bundle constructed by U.~Vetter \cite{ve} using an explicit 
matrix of linear forms.  
\vskip5mm

\subsection{Example}(Horrocks \cite{ho2})\label{nullcorr} 
Assume that  
$G^0 = \Lambda^{\vee} (-i)$ and $G^{-1} = \Lambda^{\vee} (-j-1)$, for some 
$0 < i < j < n$. The differential $d^{-1}:G^{-1} \to G^0$ is defined by a 
linear function $\Lambda^{j-i+1}V^\ast \to k$, i.e., by an element 
$\omega \in \Lambda^{j-i+1}V$. 
Let $-\cdot -:\Lambda^{\vee} \times \Lambda \to \Lambda^{\vee}$ denote the 
bilinear map defining the structure of right $\Lambda$--module on 
$\Lambda^{\vee} $, i.e., the contraction map. Then $d^{-1}_0:G^{-1}_0 \to 
G^0_0$ is just $-\cdot \omega:\Lambda^{j+1}V^\ast \to \Lambda^iV^\ast $, 
which can be identified (up to sign) with 
$-\wedge \omega:\Lambda^{n-j}V \to \Lambda^{n-i+1}V$.  

Let us give, following Horrocks \cite{ho2}, some examples of 
elements $\omega $ for which $d^{-1}_0$ is injective. Let 
$e_0,\dots ,e_n$ be a $k$--basis of $V$ and $X_0,\dots ,X_n$ the dual 
basis of $V^\ast $. Assume $n$ is odd, $n = 2m-1$, and consider the 
elements 
\[
\alpha = \mathop{\sum}\limits _{i=0}^{m-1} e_i \wedge e_{m+i} 
\in \Lambda^2V 
\quad \text{and}\quad 
\beta = \mathop{\sum}\limits _{i=0}^{m-1} X_i \wedge X_{m+i}  
\in \Lambda^2V^\ast \,.
\]
Let $\alpha^{(i)} = \frac{1}{i!}\alpha^{\wedge i} \in \Lambda^{2i}V$ be 
the $i$th divided power of $\alpha$, for $0 < i < m$. Then 
$-\cdot \alpha^{(m-i)}:\Lambda^{2m-i}V^\ast \to \Lambda^iV^\ast $ can 
be identified with $\beta^{(i)}\cdot -:\Lambda^iV \to \Lambda^iV^\ast $ 
via $\Lambda^{2m-i}V^\ast \cong \Lambda^iV$. For $i = 1$, 
$\varphi := \beta \cdot -:V \to V^\ast $ is an isomorphism. But, for 
$i > 1$, it is not so easy to show that $\beta^{(i)}\cdot -$ is an 
isomorphism (in fact, one needs some assumptions on $\text{char}k$). 
The point is that $\beta^{(i)}\cdot -$ is not equal to 
$\Lambda^i\varphi $ times a non-zero constant. 
Assume, for simplicity, that $\text{char}k = 0$. Consider the 
endomorphisms $A = -\cdot \alpha$ and $B = \beta \wedge -$ of the 
$k$--vector space $\Lambda^{\vee} $ and let $C = [A,B]$.   Then one can 
check that, for $\eta \in (\Lambda^{\vee} )_{-p} = \Lambda^pV^\ast $, 
$C\eta = (m-p)\eta $. It follows that $[C,A] = 2A$ and $[C,B] = -2B$, 
hence $A,\  B,\  C$ define a representation of Lie algebras 
$\rho : sl_2 \to gl(\Lambda^{\vee} )$ (see, for example, \cite{gh}, pp. 
118--120). One deduces that $A^i : \Lambda^{m+i}V^\ast \to 
\Lambda^{m-i}V^\ast $ is an isomorphism for $0\le i\le m$. 

We shall denote, for later use, by $\ke_i$, $0 < i < m$, the bundle 
obtained from the above $G^\bullet $ as in \ref{rkdec} 
for $n = 2m-1$, $j = 2m-i-1$ 
and $\omega = \alpha^{(m-i)}$. The bundle  $\ke_1$ is the so--called 
null--correlation bundle and has rank $n-1$ on $\P_n$, $n$ odd.       
\vskip5mm

\subsection{Example}(Famous vector bundles)\label{famous}
The Horrocks--Mumford rank 2 vector bundle on $\P_4$, \cite{hom},
can be obtained by the method from \ref{constr} and \ref{monoepi}  
starting with $G^\bullet = k^2 \otimes \Lambda^{\vee} (-2)$. 
The Sasakura rank 3 vector bundle on $\P_4$, \cite{ads}, can be also 
obtained in this way starting with a certain HT--complex $G^\bullet $ 
with $G^0 = \Lambda^{\vee} (-1) \oplus 
\Lambda^{\vee} (-2)$ and $G^{-1} = \Lambda^{\vee} (-3)$. 

More complicated, from the point of view of these general considerations, 
is the construction of the Horrocks rank 3 vector bundle on $\P_5$, 
\cite{ho3}. Assume $n = 5$, and consider the vector bundle $\ke_2$ 
defined at the end of \ref{nullcorr} (which exists, in fact, if 
$\text{char}k \not= 2$).   Since $d^{-1}_0 : G^{-1}_0 \to G^0_0$ is an 
isomorphism, the morphism $\delta^{\prime} : \kg^{-1} \to \kg^0$ from 
the proof of \ref{rkdec} is the composite morphism 
\[
G^0_{-1} \otimes \ko(-1) \lra G^0_0 \otimes \ko 
\overset{\approx}{\lla} G^{-1}_0 \otimes \ko \lra G^{-1}_1 \otimes \ko(1)\,. 
\]
One has $G^0_{-1} = \Lambda^3V^\ast = G^{-1}_1$, a decomposition 
$\Lambda^3V^\ast = B^0_{-1} \oplus Z^{-1}_1$ and $P := Z^{-1}_1$ is the set 
of primitive elements $\eta $ of $\Lambda^3V^\ast $, i.e., those with 
$A\eta = 0$. If $\eta \in P$ and $v \in V$, then $\eta \cdot v$ is a 
primitive element of $\Lambda^2V^\ast $ hence $AB(\eta \cdot v) = 
\eta \cdot v$. This means that $d^{-1}_0(\beta \wedge (\eta \cdot v)) = 
\eta \cdot v$. It follows that the differential $\kl^{-1} \to \kl^0$ of 
the complex $\kl^\bullet $ from \ref{minim} can be identified with the 
morphism $P \otimes \ko(-1) \to P \otimes \ko(1)$ whose reduced fibre 
at $[v] \in \P_5$ is $P \to P$, 
$\eta \mapsto (\beta \cdot v) \wedge (\eta \cdot v)$. 

Let, now, $\mu : \ko(-1) \to \ke_2$ be the morphism defined by the element 
$X_0 \wedge X_1 \wedge X_2 + X_3 \wedge X_4 \wedge X_5 \in P$ and let 
$\pi : \ke_2 \to \ko(1)$ be the morphism defined by  
the linear function $P \to k$ induced by 
$e_0 \wedge e_1 \wedge e_2 + e_3 \wedge e_4 \wedge e_5 \in \Lambda^3V$. 
Using the above description of the differential $\kl^{-1} \to \kl^0$, one 
cheks that $\pi \circ \mu = 0$. Moreover, one can check, using 
\ref{monoepi},(a),(c), that $\mu$ is a monomorphism of vector bundles 
and $\pi$ an epimorphism. Consequently, $\mu$ and $\pi$ define a monad 
\[
0 \lra \ko(-1) \lra \ke_2 \lra \ko(1) \lra 0
\]
whose cohomology is a rank 3 vector bundle on $\P_5$.  


\vskip1cm

\begin{thebibliography}{20}
\bibitem{ads} \textbf{H.~Abo, W.~Decker, N.~Sasakura}, An elliptic conic 
  bundle in $\P^4$ arising from a stable rank-3 vector bundle, 
  Math. Z. 229 (1998), 725--741 

\bibitem{bgg} \textbf{I.N.~Bernstein, I.M.~Gel'fand, S.I.~Gel'fand}, Algebraic
  bundles over $\P^n$ and problems of linear algebra, Funktsional'nyi Analyz i
  Ego Prilozheniya 12, No.\ 3 (1978), 66--67

\bibitem{co} {\bf I.~Coand\u{a}}, On the Bernstein-Gel'fand-Gel'fand 
correspondence and a result of Eisenbud, Fl{\o}ystad, and Schreyer, 
preprint math.AG/0206264, to appear in Journ. Math. Kyoto Univ.

\bibitem{efs} {\bf D.~Eisenbud, G.~Fl{\o}ystad, F.-O.~Schreyer}, Sheaf 
cohomology and free resolutions over exterior algebras, 
preprint math.AG/0104203, to appear in Transactions of the AMS

\bibitem{fl1} {\bf G.~Fl{\o}ystad}, Describing coherent sheaves on projective 
spaces via Koszul duality, preprint math.AG/0012263

\bibitem{fl2} {\bf G.~Fl{\o}ystad}, Koszul duality and equivalences of 
categories, preprint math.RA/0012264

\bibitem{ge} {\bf S.I.~Gel'fand}, Sheaves on $\P^n$ and problems of linear
  algebra, appendix to the Russian edition of the book Ch.~Okonek,
  M.~Schneider, H.~Spindler: Vector Bundles on Complex Projective Spaces,
  Birkh\"auser 1980, Mir, Moscow 1984

\bibitem{gh} {\bf Ph.~Griffiths, J.~Harris}, Principles of algebraic 
  geometry, John Wiley \& Sons, Inc., New York 1978

\bibitem{ho1} {\bf G.~Horrocks}, Vector bundles on the punctured spectrum of
  a local ring, Proc.\ London Math.\ Soc.\ 14 (1964), 689--713

\bibitem{ho2} {\bf G.~Horrocks}, Construction of bundles on $\P^n$, In:
  A.~Douady and J.-L.~Verdier (eds.): Les \'{e}quations de Yang--Mills,
  S\'{e}minaire E.N.S.\ (1977--1978), Ast\'{e}risque 71--72, Soc.\ Math.\ de
  France (1980), 197--203

\bibitem{ho3} {\bf G.~Horrocks}, Examples of rank three vector bundles 
  on five-dimensional projective space, J. London Math. Soc. 18 (1978), 
  15--27 

\bibitem{hom} {\bf G.~Horrocks, D.~Mumford}, A rank 2 vector bundle on 
  $\P^4$ with 15,000 symmetries, Topology 12 (1973), 63--81

\bibitem{ta} \textbf{H.~Tango}, An example of indecomposable vector bundle of
  rank $n-1$ on $\P^n$, J.\ Math.\ Kyoto Univ.\ 16, No.\ 1 (1976), 201--207

\bibitem{tr1} \textbf{G.~Trautmann}, Darstellung von Vektorraumb\"undeln
  \"uber ${\mathbb C}^n\smallsetminus\{0\}$, 
Archiv der Mathematik (Basel) 24 (1973),
  303--313 

\bibitem{tr2} {\bf G.~Trautmann}, Moduli of Vector Bundles on 
$\P_n(\mathbb C)$, Math.Ann. 237(1978), 167-186

\bibitem{ve} \textbf{U.~Vetter}, Zu einem Satz von G.~Trautmann \"uber den
  Rang gewisser koh\"arenter analytischer Moduln, Archiv der Mathematik
  (Basel) 24 (1973), 158--161
\end{thebibliography}
\end{document}